\documentclass{article}
\usepackage{graphicx} 
\usepackage{amsmath}
\usepackage[portrait, margin=3cm]{geometry}

\usepackage[square,numbers]{natbib}
\newcommand{\bh}{\ensuremath{\mathbf{h}}}
\newcommand{\bi}{\ensuremath{\mathbf{i}}}
\newcommand{\bj}{\ensuremath{\mathbf{j}}}
\newcommand{\bk}{\ensuremath{\mathbf{k}}}
\newcommand{\beem}{\ensuremath{\mathbf{m}}}

\newcommand{\bv}{\ensuremath{\mathbf{v}}}
\newcommand{\bx}{\ensuremath{\mathbf{x}}}

\newcommand{\bA}{\ensuremath{\mathbf{A}}}

\newcommand{\bP}{\ensuremath{\mathbf{P}}}

\newcommand{\bX}{\ensuremath{\mathbf{X}}}

\newcommand{\half}{\ensuremath{\frac{1}{2}}}

\newcommand{\etal}{{\it{et al}}.}

\usepackage{bm}

\newcommand{\bdu}{\bm{du}}
\newcommand{\bdv}{\bm{dv}}

\newcommand{\bxi}{\bm{\xi}}

\newcommand{\bdpsi}{\bm{d\psi}}

\newcommand{\bdA}{\bm{dA}}
\newcommand{\bdalpha}{\bm{d\alpha}}

\title{A Complete Graphic Statics for Rigid-Jointed 3D Frames. \\Part 1: Legendre Transforms for Moments}
\author{Allan McRobie, fam20@cam.ac.uk \\
Dept. of Engineering, Cambridge University, CB2 1PZ\\
https://orcid.org/0000-0002-6610-5927}

\begin{document}

\maketitle

\section*{\begin{center} Abstract \end{center}}
We extend graphic statics to describe the forces and moments in any 3D rigid-jointed frame structure. For a structure whose members are bars, graphic statics relates the {\it form diagram} (the geometrical layout of the bars) to a reciprocal {\it force diagram} which represents the forces in those bars. For 3D structures, one version of graphic statics -  {\it Rankine reciprocals} -  represents bar forces by areas of polygons perpendicular to those bars. Unfortunately, that description is incomplete: states of self-stress can exist which cannot be so represented.
 In this paper, Rankine reciprocals are generalised to provide a complete description. Not only can any state of axial self-stress be described, but so can any state of self-stress involving axial and shear forces coexistent with bending and torsional moments.
 This is achieved using a discrete version of  
 Maxwell's {\it Diagram of Stress} which maps the body space containing the structure into the stress space containing the force diagram. This mapping is a Legendre transform, defined via a stress function and its gradients. The description resulting here is applicable to any state of self-stress in any 3D bar structure whose joints may have any degree of fixity. By inclusion, the description also applies to any 2D truss or frame.

 Using elementary homology theory, a structural frame is decomposed into a set of loops, with states of self-stress being represented by dual loops in the stress space. Loops need not be plane. At any point on the structure the six components of stress resultant (axial and two shear force components, with torsional and two bending moment components) are represented by the oriented areas of the dual loops projected onto the six basis bivector planes in the 4D stress space.
These loops are the generalisation of the force polygons of Rankine’s reciprocal polyhedra. This description is complete, in that any self-stress of any frame can be represented. 

It is shown how this formulation in terms of a Legendre transform agrees with an earlier description, the Corsican sum, which was a Minkowski sum of 3D form and Rankine force diagrams. 

Finally, this paper also describes a new geometrical object whose projections encode the internal moments. This object is a hybrid of form and force: it plots the {\it original} stress function at the {\it dual} coordinates. This object allows the internal bending and torsional moments at any bar cross-section to be readily separated from the moments about the origin associated with the axial and shear forces. 

Part 1 of this paper, here, focuses on the Legendre transform and its role in representing moments. Subsequent parts will cover: using the CW-complexes of algebraic topology and homology to generalise the more traditional polyhedral form and force diagrams; the extension to describe the graphic kinematics of 3D frames;  the use of exterior algebra to reformulate the Principle of Virtual Work.

\section{Introduction}
Traditionally, graphic statics describes axial forces in 2D pin-jointed trusses. Central to this are Maxwell's reciprocal diagrams~\cite{maxwell1864a}, with line lengths in a force diagram representing force magnitudes in bars of the form diagram. Maxwell~\cite{maxwell1870a} generalised this to 3D trusses, with forces represented by polygonal areas. The resulting 3D force diagrams are sometimes called Rankine reciprocals, in recognition of Rankine's short paper on the subject~\cite{rankine1864}. Here, a further generalisation of Rankine reciprocals uses oriented loops, thereby allowing moments in frames to be represented.

This is the culmination of a sequence of publications aimed at generalising graphic statics. McRobie~\cite{McRobieRSOS1} considered 3D frames, with a construction called the Corsican sum involving a stress function and oriented loops in 4D space to describe forces and moments. 
McRobie~\etal~\cite{McRobieRSOS2} looked at associated graphical representations for structural kinematics. Various subsections of the present paper have been presented at conferences of the International Association of Shell and Spatial Structures (IASS). Specifically, the use of homological techniques to describe structural equilibrium in terms of CW-complexes was presented at IASS Boston 2018~\cite{McRobieBoston}, the extensions to cohomology and kinematics were presented at IASS Barcelona 2019~\cite{McRobieBarcelona} and  
the description of form-force duality in terms of Legendre transforms was presented at IASS Zurich 2024~\cite{McRobieZurich}. This paper builds on these to provide a single coherent description. Part 1, here, focuses on the Legendre transform and its use in representing moments. 

Other initiatives which aim to create generalised versions of graphic statics have been put forward by Karpenkov~\etal~\cite{Karpenkov2023} and Baranyai~\cite{Baranyai2024}.
The former uses a loop formalism, but is founded on homotopy rather than homology, and focuses on states of stress involving only axial force. It does not include moments or kinematics. Both use the language of exterior algebra, with differential forms defined over surfaces. 

The inclusion of moments is a strength of this paper. It can be argued that the most efficient way to carry loads is by purely axial action~\cite{BeyondBending}. Nevertheless many structures involve moments: beams are common structural elements and a theory of structures should be capable of describing them. Moreover, by including moments a theory can be built capable of describing any state of self-stress in a frame, even those containing only axial forces. Many earlier theories required stress functions to be continuous and have flat faces.  By lifting such restrictions, many of the short-comings of the earlier theories are overcome. Finally,  some trusses admit states of axial self-stress only when they are arranged in specific geometric configurations. Any deviation, even minutely, from such geometries can mean that no self-stress is possible. A theory which includes moments avoids this on-off binary: if the geometry is not quite perfect for a state of purely axial stress, then a nearby state of self-stress can be represented which is predominantly axial, but has small moments.

The following sections present the full theory. Section 2 shows how Maxwell's Diagram of Stress is a Legendre transform which defines a dual stress function which contains information about the forces. This description can represent both continuous and discrete stress fields in 2D and 3D. 
Section 3 shows that, as well as forces, Maxwell's Diagram of Stress also contains information about moments: both forces and moments are represented by the projected areas of oriented loops in a 4D stress space.  Section 4 
relates this Legendre description to the Corsican sum construction of McRobie~\cite{McRobieRSOS1}. Section 5 provides further details on the moments,  showing how projections of the original stress function allow the bending and torsional moments at a cross-section to be distilled from the overall moment defined by the dual stress function. 
Section 6 provides an example, a simple structure consisting of only a single loop. Most of the results of this paper have been presented previously~\cite{MaxwellReader, McRobieRSOS1, McRobieBoston, McRobieBarcelona, McRobieRSOS2, McRobieZurich}. The main novelty here is Section 5 on the representation of internal moments,

Part 2 of this series of papers will describe the loop formalism that will be developed using homology theory. Subsequent Parts
will extend the loop description to represent structural kinematics, and a novel description of Virtual Work will be presented. Virtual Work is a cornerstone of structural analysis, and in the description developed here, it will be recognised as a top form in 4D space.

\section{Maxwell's Diagram of Stress as a Legendre Transform}
In 1870, Maxwell~\cite{maxwell1870a} described a mapping between points in a body and points in a stress space. The mapping was defined such that the force acting on any surface in the body was equal to the force acting on the image of that surface in the stress space, even though the stress space contained only a simple, uniform isotropic pressure field. This was achieved by defining a stress function $f(x,y,z)$ whose gradient at the point $\bx = (x,y,z)$ defines a reciprocal point $\bxi = (\xi, \eta, \zeta)$ in the stress space. A dual stress function $\phi$ was defined at $\bxi$ by setting $\phi(\xi,\eta, \zeta)$ equal to the negative of the intercept of the tangent plane of $f(x,y,z)$ at $\bx$ with the $f$ axis. 
That is, for any differentiable function $f(x,y,z)$, there is a dual function $\phi(\xi,\eta, \zeta)$ defined by
\begin{equation}
\phi  =  \xi x + \eta y +  \zeta z - f    \quad
\mathrm{where} \quad \xi  =  \frac{\partial f}{\partial x}, \quad 
\eta = \frac{\partial f}{\partial y}, \quad \zeta = \frac{\partial f}{\partial z} \ .
\label{eqn1}
\end{equation}
This is a Legendre transform, and may be written more symmetrically as
\begin{equation}
    f + \phi   =  \bxi.\bx \quad
    \mathrm{with} \quad  \bxi = \mathrm{grad} \  f  \quad  \mathrm{and} \quad  \bx = \mathrm{grad} \ \phi \ .
\end{equation}

\begin{figure*}
\begin{center}
\includegraphics[width=0.8\textwidth]{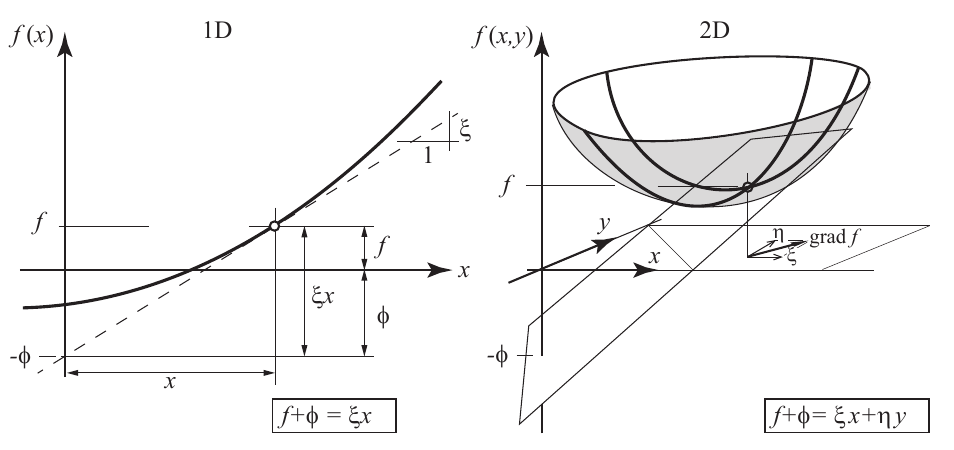}
\end{center}
\caption{The Legendre transform in 1D and 2D. The point $f(\bx)$ maps to the dual point $\phi(\bxi)$, where $\phi$ is the negative intercept of the plane tangent to $f$ at $\bx$ and $\bxi = \mathrm{grad} f$ is its gradient.}
\label{LegendreSimple}
\end{figure*}

Equation~\ref{eqn1} is for the case of stresses in 3D. The 1D and 2D cases are the obvious simplifications of this, and are illustrated in Fig.~\ref{LegendreSimple}.

In the 2D case, the stress function $f$ may be identified as the familiar Airy stress function~\cite{airy1863}, whose second derivatives define a set of stresses in equilibrium. In 3D, the stress function may be similarly identified as that associated with Rankine polyhedra, with cofactors of the matrix of second derivatives of the stress function giving a set of stresses in 3D equilibrium. 

\begin{figure*}
\begin{center}
\includegraphics[width=\textwidth]{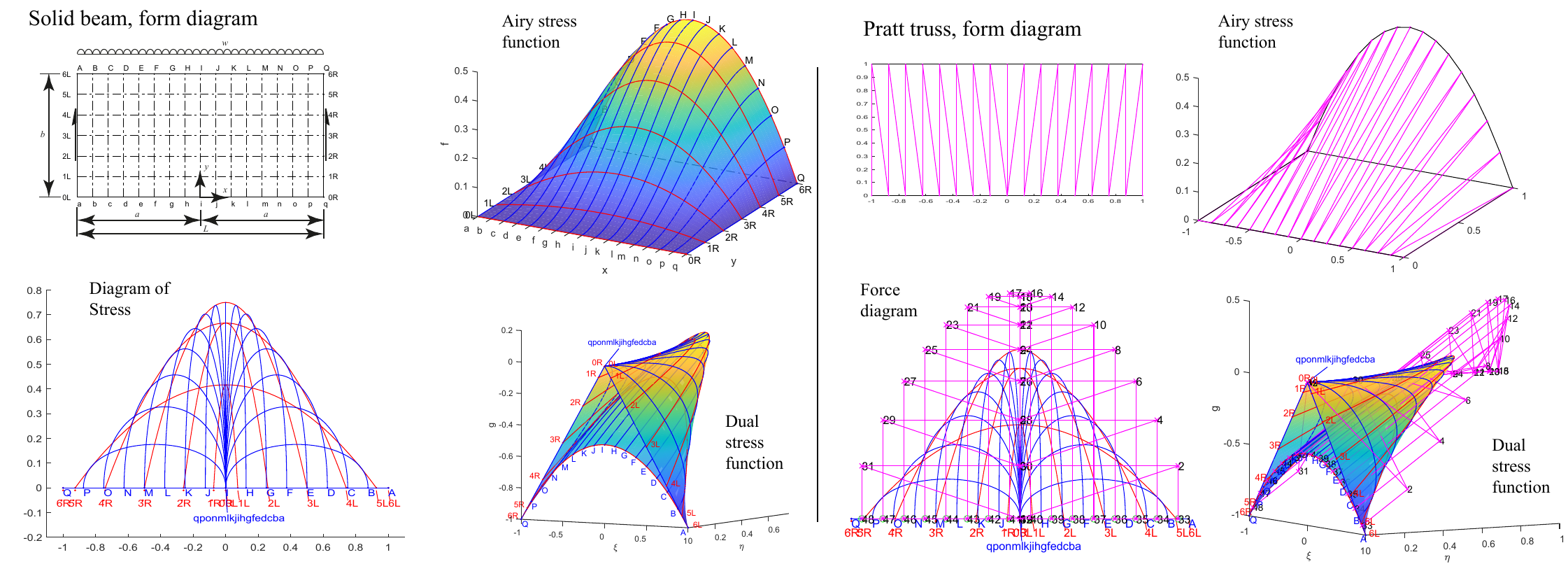}
\end{center}
\caption{A solid beam (left) and a Pratt truss (right). Each group of four diagrams shows the structure, the Airy stress function, the Diagram of Stress or force diagram, and the dual stress function. In the last two subfigures for the truss (the force diagram and the dual stress function), the results for the solid beam have been superposed for the purposes of comparison. The figure illustrates how the Legendre transform provides a continuous version of graphic statics.}
\label{BeamPratt}
\end{figure*}

Graphic statics, in both its familiar 2D form and in its 3D Rankine form, is thus a piecewise-linear manifestation of the Legendre transform. Equivalently, the Legendre transform provides a continuous version of graphic statics. 
For example, Fig.~\ref{BeamPratt} shows the Diagrams of Stress for a uniformly-loaded 2D beam and for a Pratt truss carrying the same load. 
It can be seen that the dual stress function $\phi$ is similar in the two cases. Further examples are given in \cite{MaxwellReader, McRobieZurich}. 

The relationship between the Diagram of Stress, which uses only first derivatives of the stress function, and the more usual interpretation 
 - of stresses defined via second derivatives - is described in greater detail in McRobie~\etal~\cite{McRobieZurich}.

\section{Representation of Moments}

In the following, we stay with the first derivatives of the stress function, and look beyond the stresses and forces to the moments. This additional information about moments was not mentioned by Maxwell and yet it is naturally contained within the Legendre transform that defines the Diagram of Stress. Most of the results of this section have been presented previously in ~\cite{MaxwellReader, McRobieRSOS1, McRobieBoston, McRobieBarcelona, McRobieRSOS2, McRobieZurich}.

\begin{figure*}
\begin{center}
\includegraphics[width=0.8\textwidth]{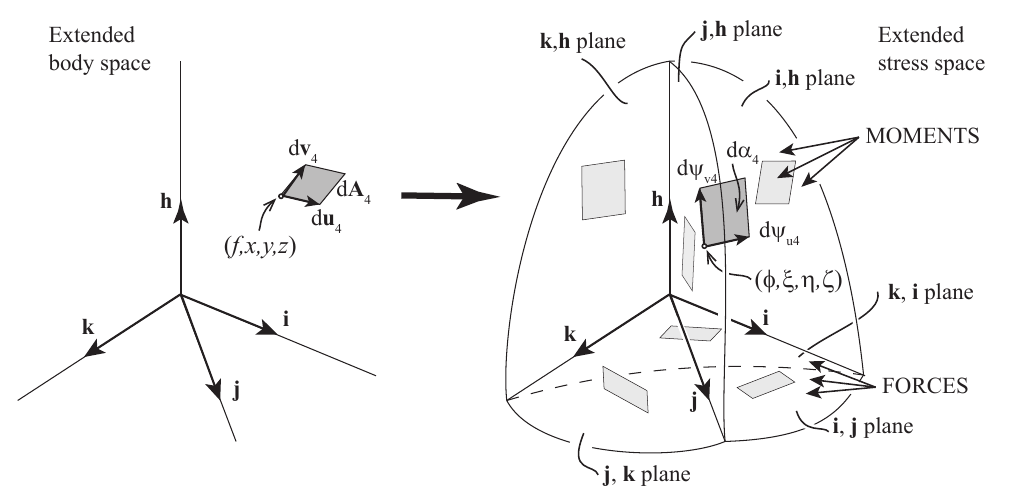}
\end{center}
\caption{A 2D element $\bdA_4$ being mapped to its image $\bdalpha_4$ in the $4\rm{D}$ setting. The six independent projections of $\bdalpha_4$ are shown schematically.}
\label{extended4D}
\end{figure*} 

Consider the $4D$ mapping from the body space point  
$(f,x,y,z)$ to the stress space point $(\phi,\xi,\eta,\zeta)$, as defined by Equation~ \ref{eqn1}. In each space we use orthonormal unit basis vectors $\bh, \bi, \bj, \bk$ in the directions $(f,x,y,z)$ of body space and $(\phi,\xi,\eta,\zeta)$ of stress space. Consider also two small vectors $\bdu_4$ and $\bdv_4$ based at the point $(f,x,y,z)$ (see Figure~\ref{extended4D})
\begin{equation*}
\bdu_4  =  df_u \bh + dx_u \bi + dy_u  \bj + dz_u \bk \quad \mathrm{and} \quad
\bdv_4  =  df_v \bh + dx_v \bi + dy_v  \bj + dz_v \bk 
\end{equation*}
 The {\it wedge product} of these two vectors defines the elemental oriented area~$\bdA_4 =  \bdu_4 \wedge \bdv_4$ of interest. 
 This small area is mapped from the 4D body space $(f,x,y,z)$  to the 4D stress space $(\phi,\xi,\eta,\zeta)$, with the images of $\bdu_4$, $\bdv_4$ and 
$\bdA_4$ being $\bdpsi_{u4}$, $\bdpsi_{v4}$ and 
$\bdalpha_4$ respectively. 
The image vectors are 
\begin{equation*}
\bdpsi_{u4}  =  d\phi_u \bh + d\xi_u \bi + d\eta_u  \bj + d\zeta_u \bk \quad \mathrm{and} \quad 
\bdpsi_{v4}  =  d\phi_v \bh + d\xi_v \bi + d\eta_v  \bj + d\zeta_v \bk 
\end{equation*}
where the magnitude of the component of $\bdpsi_{u4}$ in the $\bh$ direction is given by 
\begin{eqnarray*}
d\phi_u 
    & = & \phi(\xi+d\xi_u, \eta+d\eta_u,\zeta+d\zeta_u) - \phi(\xi, \eta, \zeta) \\
    & = & \left[(x + dx_u)(\xi + d\xi_u) +  (y + dy_u)(\eta + d\eta_u)  
     +  (z + dz_u)(\zeta + d\zeta_u)  - (f+ df_u) \right] \\ 
    & & \quad \quad \quad \quad - \left[ x\xi + y\eta + z\zeta - f \right] \\
    & \approx &  ( \xi dx_u + x d\xi_u) + ( \eta dy_u + y d\eta_u) + ( \zeta dz_u + z d\zeta_u)   
      - \left( \frac{\partial f}{\partial x} d x_u + \frac{\partial f}{\partial y} d y_u  + \frac{\partial f}{\partial z} d z_u    \right)  \\ 
   & = &  ( \xi dx_u + x d\xi_u) + ( \eta dy_u + y d\eta_u) + ( \zeta dz_u + z d\zeta_u)   
   \ - \left[ \xi d x_u + \eta d y_u  + \zeta d z_u    \right] \\
   & = &   x d\xi_u + y d\eta_u + z d\zeta_u
\end{eqnarray*}
and the second component is given by
\begin{equation*}
d\xi_u \quad   =  \quad \frac{\partial\xi}{\partial f} df_u 
       + \frac{\partial\xi}{\partial x} dx_u +
         \frac{\partial\xi}{\partial y} dy_u +
         \frac{\partial\xi}{\partial z} dz_u  \ \
     =  \quad \frac{\partial\xi}{\partial x} dx_u +
      \frac{\partial\xi}{\partial y} dy_u +
      \frac{\partial\xi}{\partial z} dz_u 
\end{equation*}
because, in the first term, 
$\partial\xi /\partial f =  (\partial/\partial f) ( \partial f/\partial x) = 0$.
The third and fourth components of $\bdpsi_{u4}$ are similar, but with $\xi$ replaced by $\eta$ and $\zeta$ respectively. The components of $\bdpsi_{v4}$ are similar, with $u$ replaced by $v$. 

The oriented area of the image element $\bdalpha_4$ is thus given by 
\begin{eqnarray*}
\bdalpha_4 & =  & \bdpsi_{u4} \wedge \bdpsi_{v4} \quad
       = \ \left(  d\phi_u \bh + d\xi_u \bi + d\eta_u  \bj + d\zeta_u \bk \right)  
      \wedge  \left(  d\phi_v \bh + d\xi_v \bi + d\eta_v  \bj + d\zeta_v \bk \right)\\
     & = & \quad (d\eta_u d\zeta_v - d\eta_v d\zeta_u) \bj \wedge \bk  
       \ \  +  \ \ \ (d\zeta_u d\xi_v - d\zeta_v d\xi_u) \bk \wedge \bi 
         \ \  + \ \ \ \ (d\xi_u d\eta_v - d\xi_v d\eta_u) \bi \wedge \bj \\
      & &  +   (d\phi_v d\xi_u - d\phi_u d\xi_v) \bi \wedge \bh  
        \ \   +       (d\phi_v d\eta_u - d\phi_u d\eta_v) \bj \wedge \bh 
      \ \   +    \   (d\phi_v d\zeta_u - d\phi_u d\zeta_v) \bk \wedge \bh 
\end{eqnarray*}

The six components of the oriented area are now clearly visible. Three are on the basis planes of the usual 3D stress space spanned by $\bi$, $\bj$ and $\bk$, where the components of the force acting on $\bdA$ in the $\bi$, $\bj$ and $\bk$ directions are given by~$p$ times the oriented areas on the $\bj \wedge \bk$, $\bk \wedge \bi$ and $\bi \wedge \bj$ planes respectively (see Fig.~\ref{extended4D}).

The components of the moment exerted about the origin by the stresses acting on $\bdA$ (the projection of $\bdA_4$ into the 3D body space) are given by the oriented areas on the other three basis planes. Specifically, the moment about the origin has components about the $\bi$, $\bj$ and $\bk$ axes given by $p$ times the components of oriented area on the $\bi \wedge \bh$, $\bj \wedge \bh$ and $\bk \wedge \bh$ planes respectively.

These results may be readily integrated from elemental areas to finite surface patches to give the following simple yet deep result: for a 2D surface patch within the 4D extended body space of $(f,x,y,z)$, not only is the total force acting on that patch given by $p$ times the oriented area of its image within the usual 3D $(\xi,\eta,\zeta) $ stress space but, further, $p$ times the oriented areas of its image when projected onto the $(\xi,\phi)$, $(\eta,\phi)$ and $(\zeta, \phi)$ planes gives the components of the total moment about the origin exerted by the stresses on that patch.

(Note that coordinates $x,y,z$ and $\xi, \eta, \zeta$ have dimensions of Length, whereas coordinates $f$ and $\phi$ have dimensions of $(\mathrm{Length})^2$.)

The fundamental objects in this analysis  are loops in 4D space. Initially, consider the discrete case, where the loops consist of  oriented polygons with straight edges. These loops need not be flat. Such a loop may be defined by its $n$ corner nodes taken in order around the loop $\lbrace \bX_1, \bX_2, \ldots  \bX_i, \ldots , \bX_n \rbrace $, with 
$\bX = (f, x,y,z) = (f, \bx)$, with the
 edges of the loop defined as the oriented straight line segments $\lbrace \bX_{i+1} - \bX_i \rbrace$. The numbering is cyclic, such that $\bX_{n+1} \equiv \bX_1 $. The oriented area $\bA$ of the loop is then given by the ``shoelace formula'' 
\begin{equation}
   \bA \equiv \half \sum_{i = 1}^n   \bX_{i+1} \wedge \bX_{i} 
   \label{eqn:loop}
\end{equation}
the result being independent of the choice of origin. The loop exists in 4D space, thus its oriented area has six bivector components. 
Although any loop may be spanned by orientable or non-orientable surfaces,
the use of the shoelace formula restricts attention to orientable surfaces. That is, Eqn.~\ref{eqn:loop} gives the oriented area of any orientable surface whose boundary is the loop. The oriented area is thus defined by the loop and does not depend upon the choice of surface spanning the loop,

As described later, any frame structure and its stress function may be decomposed into a set of loops in the 4D space of $(f,x,y,z)$. These are then mapped forward, using Maxwell's Diagram of Stress construction, to give a set of image loops in the stress space $(\phi, \xi, \eta, \zeta)$. Finally $p$ times the various projections of the image loops will define the forces and moments at points on the original loops. This is the description of equilibrium for frame structures. In later Parts of this series of papers an analogous construction will be  defined for describing the structural kinematics.

\section{Relation to Corsican Sum}
McRobie~\cite{McRobieRSOS1} provided a complete graphic statics for self-stressed frames in 3D using a construction called the ``Corsican sum''. That earlier description is now subsumed within the new Legendre transform description, as shown in the following. 

The Corsican sum is a combination of two polytopes in 4D. The first polytope represents the structure, where points have coordinates of the form $(1,x_1,x_2,x_3)$, with the additional fourth coordinate (in the first slot) having the trivial value of unity. For the other polytope, points take the form $(a_0, a_1, a_2, a_3)$.  In the notation of this paper, $(x_1, x_2, x_3) = (x,y,z)$  and $(a_1,a_2,a_3) = (\xi, \eta, \zeta)$. All that remains is to identify $a_0$.

A stress function $f(x,y,z)$ is defined over cells in the body space, and dual to each body cell is a node in the reciprocal with coordinates $(a_0, a_1, a_2, a_3)$, with the stress function gradient defining $ a_i  = \partial f/\partial x_i$ for $i = 1, 2, 3$ as usual. The Corsican sum is then akin to a Minkowski sum: it consists of attaching, at each point on the structure, a loop whose projected areas will give the stress resultant there. 

The force components (axial and two shear) are given by the three projections of the loop in the 
3D $\bi,\bj,\bk$ space as usual. Since the $(x,y,z)$ coordinates of the point of interest are fixed, only the areas of the loop defined by the reciprocal points $(a_1, a_2, a_3)$ in 3D are of interest when computing the force.

The moments are similarly given by the projected areas of the loop onto the planes $\bh \wedge \bi$, $\bh \wedge \bj$ and $\bh \wedge \bk$. Although this resembles the Legendre construction with the dual stress function, there is a key difference: with the Legendre construction, these projected areas give the TOTAL MOMENT exerted by the stresses on the section, and this includes the contribution of the force times its lever arm about the origin. In contrast, the projected areas of the Corsican construction give the INTERNAL MOMENT alone. This consists of only the bending and torsional moments, and it does not include any moment created by the offset of the force resultant.  

The Corsican sum was defined as follows~\cite{McRobieRSOS1}.  
Let the 4-polytope $P_X$ be the structure. This is simply the structure within the 3D subspace $x_0 = 1$ of the 4D setting. 
The dual polytope $P^A$ in the 4D stress space represents the state of stress in the structure. Each node $\bA^K$ of $P^A$ defines the linear stress function over the original cell $K$ of $P_X$ that is dual to node $\bA^K$.   That is, for any point $\bX$  within body cell $K$, the value of the stress function $f$ at that point is given by 
\begin{equation}
f = \bA^K.\bX = a^0 x_0 + a^1 x_1 + a^2 x_2 + a^3 x_3
\end{equation}
Since $x_0$ is unity, this is the familiar statement of a linear function over the cell.
In the notation of this paper, this is
\begin{equation}
f = a^0  + \xi x + \eta y + \zeta z
\end{equation}
from which it is clear that setting $a_0 = -\phi$ fits with the Legendre description. 

In summary, if the coordinate $a_0$ of the reciprocal polytope in the Corsican construction is identified with the negative of the dual stress function $\phi(\xi, \eta, \zeta)$ of the Legendre description, then the Corsican sum is simply a manifestation of the Legendre transform. However,  with its resemblance to a Minkowski sum, the Corsican sum contains an alternative description of the moments which may be of more use to a designer: it directly provides the information on the internal bending and torsional moments at a cross-section. The following section will show how the same information can be extracted from the Legendre description.

\section{Internal moments}
Although the Legendre description is arguably the more elegant, the way that the Corsican sum gives a direct representation of the internal bending and torsional moments may offer some advantages. A designer often wishes to know the bending moment at a point in a beam, whilst the moment about the origin of the forces there is often of lesser interest, depending as it does upon the arbitrary choice of origin. Interestingly though, a simple representation of internal moments emerges naturally from the Legendre picture, and this is consistent with the Corsican sum.

Let points on the Legendre dual object have coordinates $ \bA = (\phi, \xi, \eta,\zeta) = (\phi, \bxi$). Consider a general cut cross-section of a structural bar, and this bar forms the edge of $N$ cells. Let these body cells surround the bar in a cyclical order $1, \ldots, K, \ldots, N$ with cyclical numbering such that $N+1 \equiv 1$. The cell number, $K$, increases as cells are traversed in the order of a right-hand screw rule along the vector normal to and emerging from the cut section.
Over a typical cell, $K$, there is a stress function $f_K(x,y,z)$ which, via Legendre,  defines a dual point $A^K  = (\phi_K, \xi_K, \eta_K,\zeta_K) = (\phi_K, \bxi_K)$.  Dual to the cut section then is the reciprocal polygon connecting the corner points $\lbrace A^K \rbrace$, for $K = 1, \ldots, N$. The  resultant of the stresses acting on the cut section then consists of ($p$ times) the various projected areas of the dual loop:  
\begin{eqnarray}
    \mathrm{Force}& = & \half p  \sum_{K=1}^N \bxi^K \wedge \bxi^{K+1} \\
    \mathrm{ and \ Total \ Moment} & = & \half p \sum_{K=1}^N  (\Phi^{K+1} + \Phi^{K}) \bh \wedge(\bxi^{K+1} - \bxi^K)  \label{eqn:TotMom} 
\end{eqnarray}
The latter equation follows because the area of the trapezoid along the edge connecting reciprocal node $K$ to $K+1$ is the wedge product of the edge vector $\bxi^{K+1} - \bxi^K$ with a vector in the $\bh$ direction of length equal to the average value of $\Phi$ along that edge.

The force expands to 
\begin{eqnarray}
\mathrm{Force} & = &  \half p  \sum_{K=1}^N  (\xi^K \bi + \eta^K \bj + \zeta^K \bk) \wedge  
    (\xi^{K+1} \bi + \eta^{K+1} \bj + \zeta^{K+1} \bk) \\
 & = &  \half p \sum_{K=1}^N  \quad  [\xi\eta] \ \bi \wedge \bj  + 
                        [\eta\zeta] \ \bj \wedge \bk  + 
                        [\zeta\xi] \ \bk \wedge \bi \label{eqn:crossprod}\\
                       \mathrm{where \ the \ symbol} \ \ [\xi\eta] & \equiv & \xi^{K}\eta^{K+1} - \xi^{K+1} \eta^{K} \ , \ \ \mathrm{etc.}
\end{eqnarray}
Notice that Eqn.~\ref{eqn:crossprod} is simply the sum of vector cross products, but expressed in bivector form. 
Alternatively, it may be represented by the Hodge dual of this bivector, the force vector $\bP$  normal to it in 3D, namely:
\begin{eqnarray}
\mathrm{Force \ vector} \ \bP  =   \half p \sum_{K=1}^N  \quad   
                        [\eta\zeta] \ \bi  + 
                        [\zeta\xi] \ \bj+ [\xi\eta] \ \bk  \\
                         = \half p \sum_{K=1}^N \quad \bxi^K \times \bxi^{K+1} \\
                         = \half p
\begin{vmatrix}
\bi & \bj & \bk\\
\xi^K & \eta^K & \zeta^K \\
\ \xi^{K+1} & \eta^{K+1} & \zeta^{K+1}
\end{vmatrix}	                  
\end{eqnarray}
This illustrates that the terms of the form $[\xi\eta]$ are simply the $2 \times 2$ determinants familiar in the more usual vector representation.

The bivector representing the Total Moment may be decomposed as
\begin{equation}
\mathrm{Total \ Moment}   =   \bx \wedge \bP  +  \mathrm{Internal \ Moment}  
\end{equation}
We shall need the readily-provable result that 
\begin{eqnarray}
\bx \wedge \bP  & = &  \half \sum_{K=1}^N  
 (y[\xi\eta] - z[\zeta\xi]) \bj \wedge \bk + 
 (z[\eta\zeta] - x[\xi\eta]) \bk \wedge \bi +
 (x[\zeta\xi] - y[\eta\zeta]) \bi \wedge \bj
 \label{Eqn:bwedgeP}
\end{eqnarray}
From Eqn.~\ref{eqn:TotMom} we have
\begin{eqnarray}
  \mathrm{ Total \ Moment}   & = & \half p  \sum_{K=1}^N  (\phi^{K+1} + \phi^{K}) \bh \wedge (\bxi^{K+1} - \bxi^K) \hspace*{3cm} \\
& = & \half p \sum_{K=1}^N    [\phi \xi] \bi \wedge \bh
 + [\phi \eta] \bj\wedge\bh +
  [\phi \zeta] \bk\wedge\bh \\
& = & \half p  \sum_{K=1}^N    [\phi \bxi] \wedge \bh \label{eqn216}
\end{eqnarray}

\begin{eqnarray}
\mathrm{Now} \quad \Phi^K & = &  x \xi^K  + y \eta^K + z\zeta^K - f_K(\bx) \\
\mathrm{so} \quad [\Phi \xi ] & = & \quad  
       ( x\xi^{K+1} + y \eta^{K+1} + z\zeta^{K+1} - f_{K+1}) \xi^K \\
& &  \quad \quad  \quad - ( x\xi^{K}   + y \eta^{K} + z\zeta^{K} - f_K) \xi^{K+1} \\
& = &  y [\xi\eta] + z[\zeta \xi] - f_{K+1}\xi^K + f_K\xi^{K+1} 
\end{eqnarray}
Repeating for the other two components, and using Eqn.~\ref{Eqn:bwedgeP}, leads to 
\begin{equation}
    \mathrm{Total \ Moment} = \bx \wedge \bP + \half p \sum_{K = 1}^N [f \bxi ] \wedge \bh
\end{equation}

whence the result
\begin{equation}
    \mathrm{Internal \ Moment} \quad  = 
    \half p \sum_{K = 1}^N (f_K \bxi^{K+1} - f_{K+1} \bxi^K) \wedge \bh = \half p  \sum_{K = 1}^N  [f \bxi] \wedge \bh
    \label{Eqn:InternalM}
\end{equation}

The projection onto the $\bi \wedge \bh$ plane gives the internal moment about the $x$ axis, which is thus 
\begin{equation}
    \half p \sum_{K = 1}^N (f_K \xi^{K+1} - f_{K+1} \xi^K) = \half p \sum_{K = 1}^N [f \xi]
\end{equation}

The polygon projected onto this $\bi \wedge \bh$ plane has corner nodes defined by vectors of the form $\bv^K = \xi^K \bi +  f^K\bh$. The shoelace formula gives the projected area $\bA_{\bi \wedge \bh}$ as
\begin{equation}
 \bA_{\bi \wedge \bh}  =  \half \sum_{K = 1}^N   \bv^{K+1}\wedge \bv^{K}  =  \half \sum_{K = 1}^N (f_K \xi^{K+1} - f_{K+1} \xi^K)  \bi \wedge \bh  = \half \sum_{K = 1}^N [f\bxi]\wedge \bh
 \label{eqn224}
\end{equation}
since $\bi\wedge \bi = \bh \wedge \bh = 0$ leaves only the cross-terms. 
This result shows that, whilst the Total Moment is given by the $\bh$ plane projections of a face of the $(\phi, \bxi)$ polytope (Eqn.~\ref{eqn216}), the Internal Moment is given by $\bh$ plane projections of a polygon with corners of the form $(f, \bxi)$ (Eqn.~\ref{eqn224}).   The overall Legendre dual polytope with corners $(\phi, \bxi)$ describes the state of self-stress of the whole structure, whereas the polygons with corners of the mixed form  $(f, \bxi)$ vary as the chosen point $\bx$ of the structure is varied. That is why the Corsican sum construction defines a dual polygon at every point on the structure. The Legendre construction, in contrast, defines a dual polygon for every structural loop, leading to a single dual polytope with which to represent the stress resultants. 

In summary, as well as the dual polytopes $(f,\bx)$ and $(\phi, \bxi)$  a third object may be defined, the hybrid object $(f, \bxi)$. This plots the original stress function $f$ at the dual coordinates $\bxi$. Projections of polygons onto basis planes of the 3D space spanned by $\bi, \bj, \bk$ give the force components, as usual. Projections of polygons onto the three additional basis planes $\bi\wedge\bh$, $\bj\wedge\bh$, $\bk\wedge\bh$ of the 4D space give the components of the internal bending and torsional moments at any cut cross-section of the structure.  

It follows that a dual hybrid object $(\phi, \bx)$ can be defined which would give internal moments if the roles of structure and dual were to be swapped. We have yet to use this possibility.

An example is now given showing how the above formalism may be used to encode any state of self-stress in a structure which consists of only a single loop. A subsequent Part of this series of papers shows how any self-stressed frame may be decomposed into a set of structural loops, thus any state of self-stress may be represented by  the union of a set dual loops. Analysis of any frame is thus simply the union of a set  of examples like the single loop example that follows.

\section{Example: a rectangular frame}

\begin{figure*}[h]
\begin{center}
\includegraphics[width=0.85\textwidth]{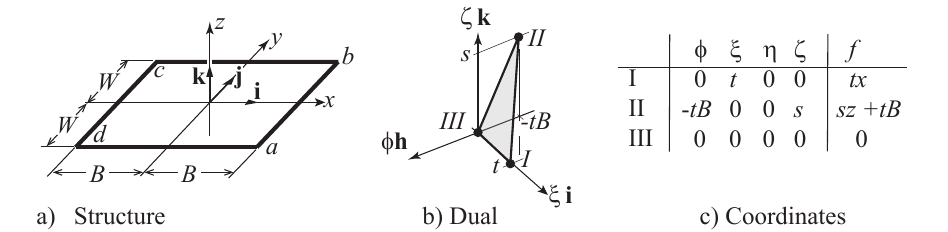}
\end{center}
\caption{ $a$) A rectangular frame. $b,c$) Dual object and coordinates. }
\label{LWsquareb}
\end{figure*}

For the rectangular frame $abcd$ of dimensions $2B$ by $2W$ in the $\bi, \bj$ directions respectively (Fig.~\ref{LWsquareb}$a$), consider 
the state of self-stress defined by a triangular loop with nodes $I-III$ in the dual space  $(\phi, \xi, \eta, \zeta)$ (Fig.~\ref{LWsquareb}$b,c$). Any loop in 4D would be admissible, and there is no particular reason for this choice, other than simplicity.  Often, when working with stress functions, a stress function $f$ is chosen, and one then sees what stresses it defines. The situation here is similar, but instead we begin by choosing the dual stress function $\phi$ (Fig.~\ref{LWsquareb}$b$). 

(Note: for brevity, the bivector $\bi \wedge \bj$ may also be written as the Clifford product $\bi\bj$.)

\begin{figure*}[ht]
\begin{center}
\includegraphics[width=0.9\textwidth]{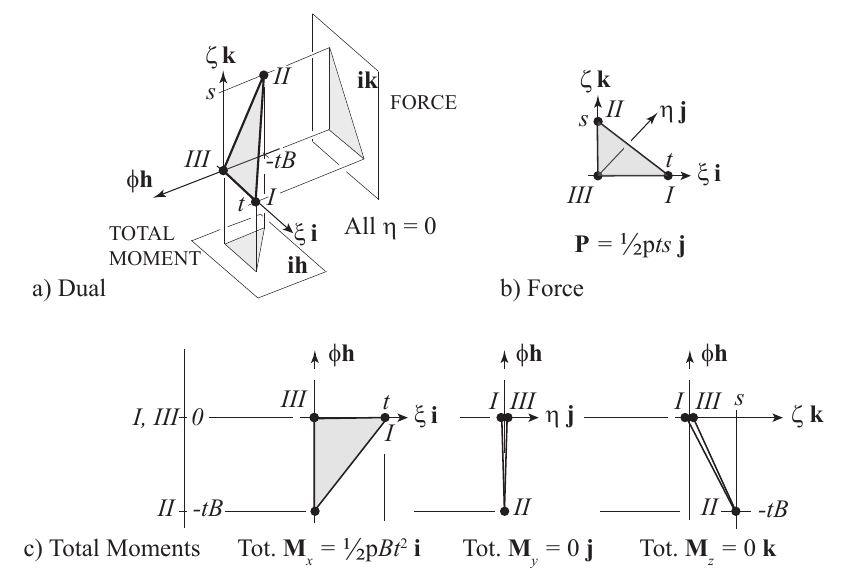}
\end{center}
\caption{ $a$) The dual object and its projections $b$) Force. $c$) Total moment. }
\label{LWsquarec}
\end{figure*}

The various projected areas of the dual object are immediately apparent (Fig.~\ref{LWsquarec}). On the $\bi\bk$ plane, the projected triangle
with sides $s$ by $t$ will give a force $\half pts$ along the $\bj$ direction. On the $\bi\bh$ plane, the triangle with sides $t$ by $tB$  will give a Total Moment  of $\half pBt^2$ about the $\bi$ axis. 

\begin{figure*}[ht]
\begin{center}
\includegraphics[width=0.9\textwidth]{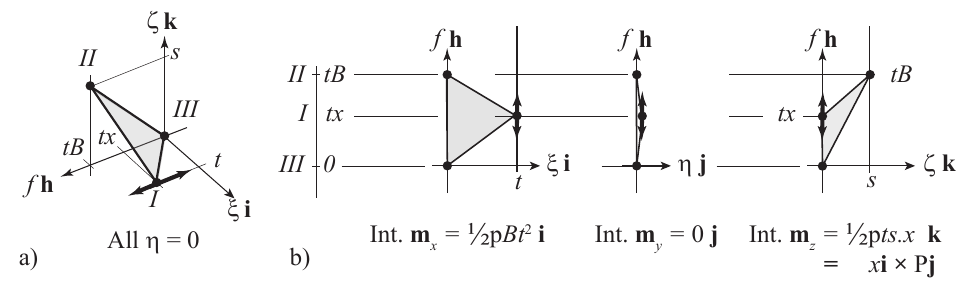}
\end{center}
\caption{ $a$) Hybrid object, plotting the original stress function $f$ at the dual coordinates $(\xi, \eta, \zeta)$.\\ $b$) The projections which give the internal moments.}
\label{Hybrid}
\end{figure*}

The dual object chosen thus represents a state of self-stress where there is a force $P = \half p t s$ in the $\bj$ direction and a Total Moment  of $\half pBt^2$ in the $\bi$ direction. This is true for all points on the structure. 

If, instead of $\phi$, the original stress function $f$ is plotted in the $\bh$ direction (Fig.~\ref{Hybrid}$a$), the projections involving the $\bh$ direction give the internal moments (Fig.~\ref{Hybrid}). 
It can be seen that although the Total Moment had no component in the $\bk$ direction, there are internal bending moments in that direction of magnitude $Px$.

\begin{figure*}[ht]
\begin{center}
\includegraphics[width=0.9\textwidth]{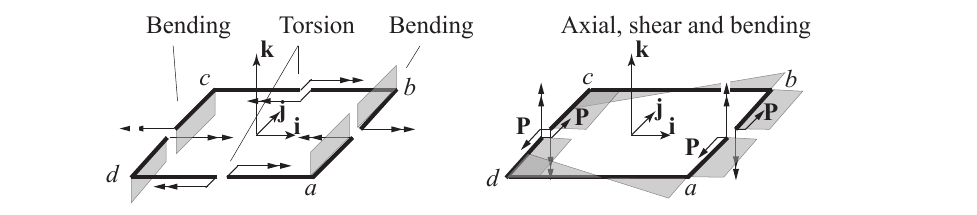}
\end{center}
\caption{Traditional representations of internal forces and moments.}
\label{SquareTrad}
\end{figure*}

Traditional representations in terms of bending moment diagrams are shown in Fig.~\ref{SquareTrad}.

The following qualitative features may be observed.
\begin{itemize}
\item The dual object is able to represent forces and moments in the fully-3D case, this being an advance over the Williams-McRobie~\cite{WilliamsMcRobie2016} use of discontinuous Airy stress functions for describing in-plane-only forces and moments in 2D frames.
\item The dual object is simple - it is a triangle - and yet it encodes a comparatively complicated state of self-stress involving axial and shear forces coexistent with torsional moments and in-plane and out-of-plane bending. The usual way of encoding this state of self-stress would be to specify two vectors, a force and a moment, at some given cut (specified by a third vector). That traditional description is also simple. The new description does not need to specify which cut is being considered, as the dual object applies to the whole loop.
\item To be fully general, the dual object would need to have at least five sides. If there are fewer sides, as here, the force and moment will be orthogonal, as here. This was proved in McRobie~\cite{McRobieRSOS1}.
\item As well as the dual diagram, a new object has been described. This is a hybrid diagram (Fig.~\ref{Hybrid}) which plots the original stress function $f$ at the dual coordinates $(\xi, \eta,\zeta)$.
\item The dual diagram, defined by plotting the dual stress function $\phi$, is fixed. It applies to all points on the structure. This follows because, by elementary equilibrium, the Internal Force $\bP$ acting across any cut section is the same at all points on the loop. Similarly the Total Moment (which includes the moment $\bx \times \bP$  associated with the Internal Force $\bP$) is also constant for all points on the loop. This follows trivially from moment equilibrium of a length of bar: there being no applied loads, the total moments at the two ends must equilibrate to zero.
\item In contrast, the geometry of the hybrid diagram depends on which point of the structure is under consideration. In this example, the variation is via the $x$ coordinate, leading to the linear variation of the bending moment $\beem_z$ along sides $da$ and $cb$ associated with the Internal Force $\bP$. 
\item The dual diagram defines a state of stress that is independent of the structure's geometry. The same dual diagram can be associated with any number of different structural geometries, as demonstrated in the following example.  
\end{itemize}

\begin{figure*}[ht]
\begin{center}
\includegraphics[width=\textwidth]{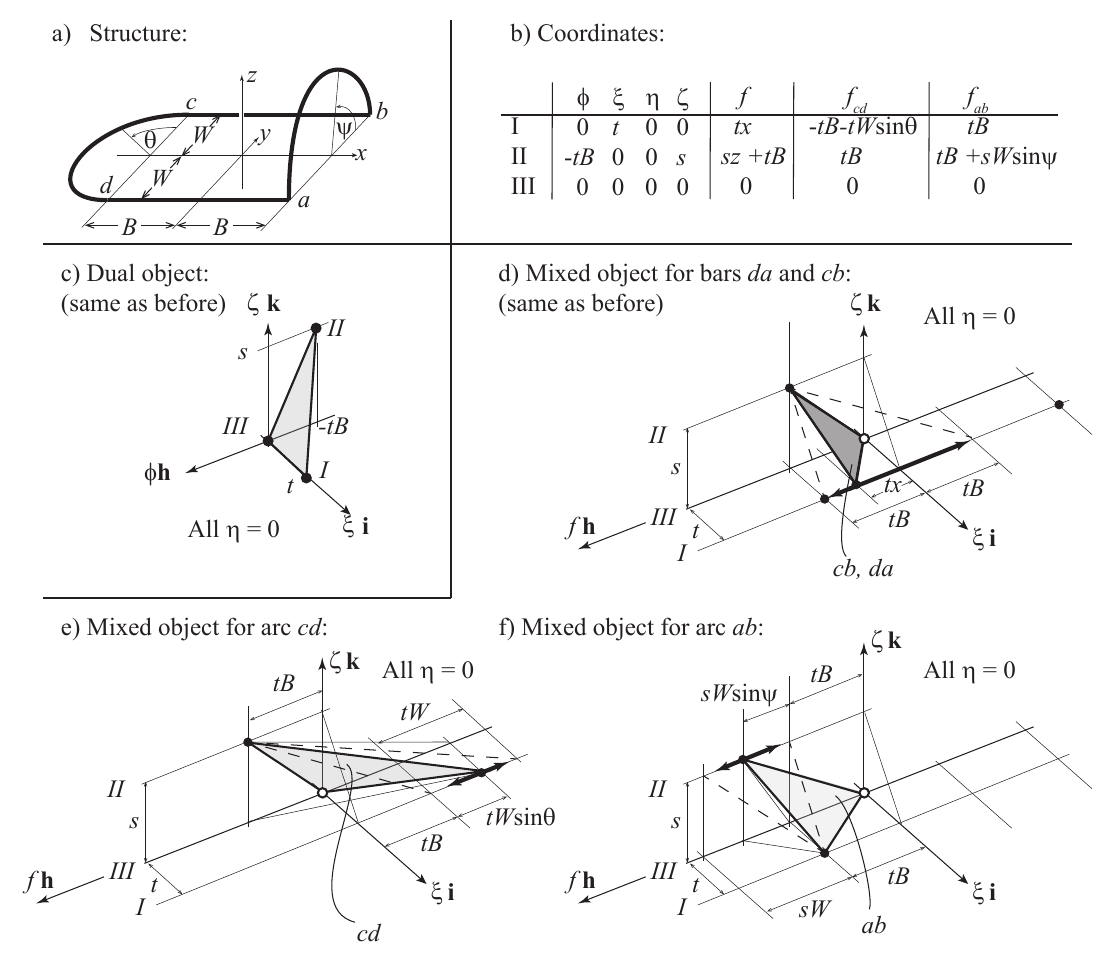}
\end{center}
\caption{ $a$) The rectangular frame with semi-circular ends. $b$) Dual coordinates. $c$)  Dual object. $d,e,f$ Mixed objects for points on various bars. }
\label{CircEnds}
\end{figure*} 

\subsection{Example with curved members}

Consider a variant of the previous example. The structure again consists of a single loop, but now two of the members are semicircular arcs,  one of which curves out of the plane (Fig.~\ref{CircEnds}$a$).  The structure thus has the same topology as previously but with different geometry. Forces and Total Moments do not vary around a single-loop structure, and these are defined by the dual loop. It follows that if the same dual object is used as previously, the Forces and Total Moments will be identical to those in the first structure. However the Internal Moments depend on the structural geometry, and so these are now different. They are nevertheless readily computed by taking projections of the mixed objects, where the stress function $f$ is plotted with respect to the dual coordinates $\xi$, $\eta$, $\zeta$ (Fig.~\ref{CircEnds}$d,e,f$). For example,  it can be read directly from Fig.~\ref{CircEnds}$f$ that for a section at an angle $\psi$ on the out-of-plane arc $ab$ the moment about the $\bk$ axis has a magnitude $\half pstB $ as previously, whereas the moment about the $\bi$ axis now has magnitude $\half pt(tB + sW \sin  \psi) $.

\section{Summary, Conclusions and Further Work}
It has been shown how Maxwell's Diagram of Stress corresponds to a Legendre transform between a stress function over the body space and a dual stress function over the stress space. It has also been shown how this description contains information about moments in structural bars. Within this description, all six components (axial and two shear forces, and torsional and two bending moments) of any stress-resultant at a section may be represented by the projected areas of a loop within the 4D stress-space. It was further demonstrated how a hybrid object may be constructed by plotting the original stress function at the dual coordinates, and that projections of such hybrid loops give direct information about the internal torsional and bending moments, decoupling these from the moments due to the lever arm about the original of the axial and shear forces. Whilst each dual loop is dual to a whole structural loop, each hybrid loop is dual to a single point on the structure.

An example was given of a self-stressed frame structure  consisting of only a single loop. In the first instance this was a
flat rectangle, but this was later varied to be a more general space curve, with some members curved.

Clearly, more complicated structures can be constructed as assemblages of such loops, with each structural loop having an associated dual loop to represent forces and moments. 
The way that any frame structure can be decomposed into loops using elementary homology theory was presented in McRobie~\cite{McRobieBoston, McRobieBarcelona, McRobieZurich}. This will be the subject of Part 2 of this sequence of papers. In essence, the 2D Maxwell and 3D Rankine theories of
graphic statics are associated, respectively, with polyhedral and polytopic stress functions, but these primitive geometrical objects have restrictions, such as requiring faces to be plane. Such restrictions mean that some structures and their self-stresses cannot be represented in those formalisms. Such limitations can be avoided by using a more general family of objects, the CW-complexes of algebraic topology~\cite{McRobieBoston, McRobieBarcelona, McRobieZurich}. These have fewer geometrical restrictions, such that any self-stress in any frame may then be represented.  In particular, bars can be general space curves and faces need not be flat.

Throughout this paper, the only geometric objects considered are loops and points on loops. Dual to each structural loop is a force loop.  There was no mention in this paper of higher dimensional objects such as faces, cells or polytopes. Such objects - which will be part of the CW-complexes -  are the generalised ``liftings'' of the form and force loop assemblages. Such liftings are not necessary for a complete theory of 3D frames - the loops alone suffice. Nevertheless, the higher dimensional objects can be useful for understanding how forces are organised, and such aspects will be considered in Part 2.



\bibliography{McRobieBib}
\end{document}